\documentclass[conference]{IEEEtran}
\IEEEoverridecommandlockouts
\usepackage{cite}
\usepackage{amsmath,amssymb,amsfonts}
\usepackage{algorithmic}
\usepackage{graphicx}
\usepackage{textcomp}
\usepackage{xcolor}
\usepackage{multirow}

\newcommand{\bbb}{\mbox{\boldmath $b$}}

\newcommand{\gggg}{\mbox{\boldmath $g$}}

\newcommand{\rrr}{\mbox{\boldmath $r$}}

\newcommand{\vvv}{\mbox{\boldmath $v$}}
\newcommand{\www}{\mbox{\boldmath $w$}}
\newcommand{\xxx}{\mbox{\boldmath $x$}}
\newcommand{\yyy}{\mbox{\boldmath $y$}}
\newcommand{\zzz}{\mbox{\boldmath $z$}}
\newcommand{\AAA}{\mbox{\boldmath $A$}}

\newcommand{\DDD}{\mbox{\boldmath $D$}}

\newcommand{\HHH}{\mbox{\boldmath $H$}}

\newcommand{\LLL}{\mbox{\boldmath $L$}}
\newcommand{\MMM}{\mbox{\boldmath $M$}}

\newcommand{\UUU}{\mbox{\boldmath $U$}}

\newcommand{\ZZZ}{\mbox{\boldmath $Z$}}

\newcommand{\beq}{\begin{equation}}
\newcommand{\eeq}[1]{\label{#1} \end{equation}}
\newcommand{\beqa}{\begin{eqnarray}}
\newcommand{\eeqa}[1]{\label{#1} \end{eqnarray}}
\newcommand{\bmat}[1]{\left ( \begin{array}{#1}}
\newcommand{\emat}{\end{array} \right )}

\def\BibTeX{{\rm B\kern-.05em{\sc i\kern-.025em b}\kern-.08em
    T\kern-.1667em\lower.7ex\hbox{E}\kern-.125emX}}
\begin{document}

\title{An Integer Arithmetic-Based Sparse Linear Solver Using a GMRES Method and Iterative Refinement
\thanks{This work was supported by JSPS KAKENHI Grant Numbers 20K21782 and 19H01105}
}

\author{\IEEEauthorblockN{Takeshi Iwashita}
\IEEEauthorblockA{Information Initiative Center\\
Hokkaido University\\
Sapporo, Japan\\
Email: iwashita@iic.hokudai.ac.jp}
\and
\IEEEauthorblockN{Kengo Suzuki}
\IEEEauthorblockA{Department of Electrical Engineering\\
Hokkaido University\\
Sapporo, Japan\\
Email: kiken50627@eis.hokudai.ac.jp}
\and
\IEEEauthorblockN{Takeshi Fukaya}
\IEEEauthorblockA{Information Initiative Center\\
Hokkaido University\\
Sapporo, Japan\\
Email: fukaya@iic.hokudai.ac.jp}
}

\twocolumn[
{\Huge IEEE Copyright Notice}

\vspace{2\baselineskip}
Copyright  \copyright \, 2020 IEEE

Personal use of this material is permitted. Permission from IEEE must be obtained for all other uses, in
any current or future media, including reprinting/republishing this material for advertising or promotional
purposes, creating new collective works, for resale or redistribution to servers or lists, or reuse of any
copyrighted component of this work in other works.

\vspace{2\baselineskip}
Accepted to be published in: 2020 IEEE/ACM 11th Workshop on Latest Advances in Scalable Algorithms for Large-Scale Systems (ScalA)

\vspace{2\baselineskip}

Cite as:

\begin{quote}
T. Iwashita, K. Suzuki and T. Fukaya, ``An Integer Arithmetic-Based Sparse Linear Solver Using a GMRES Method and Iterative Refinement," 2020 IEEE/ACM 11th Workshop on Latest Advances in Scalable Algorithms for Large-Scale Systems (ScalA), GA, USA, 2020, pp. 1-8, doi: 10.1109/ScalA51936.2020.00006.
\end{quote}

\vspace{2\baselineskip}
BibTex:

\begin{quote}
@INPROCEEDINGS\{9308712,
  author=\{T. \{Iwashita\} and K. \{Suzuki\} and T. \{Fukaya\}\},
  booktitle=\{2020 IEEE/ACM 11th Workshop on Latest Advances in Scalable Algorithms for Large-Scale Systems (ScalA)\}, 
  title=\{An Integer Arithmetic-Based Sparse Linear Solver Using a GMRES Method and Iterative Refinement\}, 
  year=\{2020\},
  volume=\{\},
  number=\{\},
  pages=\{1-8\},
  doi=\{10.1109/ScalA51936.2020.00006\}\}
\end{quote}

]

\newpage
\maketitle

\begin{abstract}
In this paper, we develop a (preconditioned) GMRES solver based on integer arithmetic, and introduce an iterative refinement framework for the solver.
We describe the data format for the coefficient matrix and vectors for the solver that is based on integer or fixed-point numbers.
To avoid overflow in calculations, we introduce initial scaling and logical shifts (adjustments) of operands in arithmetic operations. 
We present the approach for operand shifts, considering the characteristics of the GMRES algorithm.
Numerical tests demonstrate that the integer arithmetic-based solver with iterative refinement has comparable solver performance in terms of convergence to the standard solver based on floating-point arithmetic.
Moreover, we show that preconditioning is important, not only for improving convergence but also reducing the risk of overflow.
\end{abstract}


\section{Introduction}
In recent years, it has become difficult to improve the performance of processors, particularly their energy efficiency.
The main reason is the decline in lithographic scaling, which threatens the well-known techno-economic model for the IT industry, that is, Moore's Law~\cite{Vetter,Shalf}. 
Thus, new computing technologies and devices based on different physics from CMOS technology are being widely investigated.
Although quantum computing is a typical example for these technologies, some technologies aim to develop an ultra low-power but high-performance computer that is operated by instructions similar to conventional computers, for example computing devices based on single-flux-quantum (SFQ) circuits~\cite{Sato,ishida}.
However, these new types of computers may support only integer arithmetic in the early stage of research and deployment, because circuits for floating-point (FP) arithmetic are more complex and power consuming than those for integer arithmetic. Accordingly, we attempt to evaluate the potential of integer arithmetic computing for scientific computing. 
Specifically, we focus on iterative methods that are widely used in various scientific simulations, and investigate an integer arithmetic-based iterative linear solver, in which only integer arithmetic is used in the main iteration loop.
%
%
%

While there is a wide variety of iterative solvers, we develop a generalized minimal residual (GMRES) solver~\cite{Saad} using integer (fixed-point number) arithmetic that is denoted by int-GMRES.
The GMRES method is a Krylov subspace method and is used as a standard solver for a linear system that has an unsymmetric coefficient matrix.
In our solver, the iterative refinement technique is used with the GMRES solver based on integer arithmetic to
obtain a solution vector with the same accuracy as the output of a standard FP arithmetic solver.
Although the technique is classical, it is useful for mixed-precision computing~\cite{Goddeke}.
In this paper, we introduce the iterative refinement framework for an integer arithmetic-based solver and present the details of the implementation of the int-GMRES solver.

In Sections II and III, we introduce some notation and problem definitions, including the initial scaling of the linear system to be solved.
In Section IV, we describe the iterative refinement framework for the solver based on integer arithmetic.
In Sections V and VI, we present the details of the implementation of (preconditioned) int-GMRES.
In Section VII, we present the numerical results. In Sections VIII and IX, we describe related works and summarize the paper.

\section{Notation}
In this paper, we discuss a linear solver in which integer arithmetic is mainly used. In the program, some variables and elements of arrays are declared as integer numbers, and they are treated as fixed point numbers in the analysis. We use Q notation for the fixed-point number. Q$d_{m}$.$d_{f}$ denotes a number with  $d_{m}$ integer and $d_{f}$ fractional bits. The word length $WL$ is $d_{m}+d_{f}+1$, because a sign bit is used. The entire word is a two's complement integer. 

In the following, we denote the $i$-th row $j$-th element of matrix $\ZZZ$ by $z_{ij}$ or $Z(i, j)$.
We denote the $i$-th element of vector $\zzz$ by $z_{i}$.
When matrices, vectors, and variables have a {\it bar}, such as $\bar{\AAA}^{(k)}$, this indicates that their elements or values are fixed-point or integer numbers and are stored using the intWL type in the program. 

\section{Problem and Initial Scaling}
In this paper, we consider the following $n$-dimensional linear system of equations:
\beq
\hat{\AAA} \hat{\xxx} = \hat{\bbb}.
\eeq{org}
The elements of $\hat{\AAA}$ and $\hat{\bbb}$  are given by FP numbers. Typically, they are double precision. 
We need to solve (\ref{org}) with sufficient accuracy; that is, the relative residual norm calculated using (double-precision) FP arithmetic must be smaller than a given tolerance. 
The final value of each element of $\hat{\xxx}$ is given by an FP number.

First, the linear system (\ref{org}) is scaled using FP arithmetic as follows:
\beq
\AAA \xxx = \bbb,
\eeq{eq1}
where $\AAA=\hat{\DDD}^{-1}\hat{\AAA}$ and $\bbb=\hat{\DDD}^{-1}\hat{\bbb}$.
When we intend to preserve a particular property of the coefficient matrix, such as symmetry, the scaled linear system can be written as follows:
\beq
\AAA=\hat{\DDD}_{1}^{-1}\hat{\AAA}\hat{\DDD}_{2}^{-1}, \ \xxx=\hat{\DDD}_{2}\hat{\xxx}, \ \bbb=\hat{\DDD}_{1}^{-1}\hat{\bbb}. 
\eeq{eq2}
In (\ref{eq1}) and (\ref{eq2}), $\hat{\DDD}$, $\hat{\DDD}_{1}$, and $\hat{\DDD}_{2}$ are diagonal matrices.
In the present analysis, the $i$-th diagonal element of $\hat{\DDD}$ is given by
\beq
\hat{d}_{i i}=\max_{j} | a_{ij}    | / 2^{\bar{\alpha}_{a}}.
\eeq{eq3}
When the linear system (\ref{eq1}) is solved mainly using integer arithmetic, the setting of  $\bar{\alpha}_{a}$ can be an important issue, and depends on the solver implementation.
Based on our preliminary tests, we suggest that $\bar{\alpha}_{a} = WL/4$, whereas a larger value can be set for a preconditioned solver.

\section{Iterative Refinement}
We use an iterative refinement technique, which is slightly adjusted 
 for iterative linear solvers based on integer arithmetic. In the technique, we refine the approximate solution vector by solving the residual equation.
We assume that we obtain sufficiently accurate solution vector by $k_{t}$ times refinements.
In each refinement, a linear system of equations is approximately solved.
Then, the solution vector (or its sufficiently accurate approximation) $\xxx$ is written as
\beq
\xxx=\tilde{\xxx}^{(1)}+\tilde{\xxx}^{(2)}+\cdots+\tilde{\xxx}^{(k_{t})}.
\eeq{solution}
In our technique, the approximate solution vector $\tilde{\xxx}^{(k)}$ for the $k$-th refinement is obtained by (approximately) solving the linear system of equations:
\beq
\bar{\AAA}^{(k)} \xxx^{(k)} = \bbb^{(k)}.
\eeq{solve}
In (\ref{solve}), each element of $\xxx^{(k)}$ and $\bbb^{(k)}$ is given by an FP number. 
 
\subsection{Setting of the Right-Hand Side and Solution Vector}
Before the $k$-th refinement process, we calculate 
\beq
\bbb'^{(k)} =  \bbb - \AAA (\sum_{l=1}^{k-1}\tilde{\xxx}^{(l)})
\eeq{b}
using FP arithmetic.
We note that $\bbb'^{(1)}=\bbb$.
Although $\tilde{\xxx}^{(k)}$ can be determined by solving $\bar{\AAA}^{(k)} \tilde{\xxx}^{(k)} = \bbb'^{(k)}$,
we solve its scaled system (\ref{solve}) considering the use of an integer arithmetic-based solver and representation range of a fixed-point number.  
%
Using FP arithmetic, we calculate the scaled vector $\bbb^{(k)} $ of $\bbb'^{(k)}$ using
\beq
\bbb^{(k)} = \frac{1}{\gamma^{(k)}} \bbb'^{(k)},
\eeq{scale} 
and
\beq
\gamma^{(k)}=\max_{i} | b'^{(k)}_{i}  |. 
\eeq{gammak}
Then, 
the vector for the $k$-th refinement $\tilde{\xxx}^{(k)}$ is written by 
\beq
\tilde{\xxx}^{(k)}=\gamma^{(k)} \xxx^{(k)}.
\eeq{update}
When the entire refinement process works, we can expect that the scaling factor $\gamma^{(k)}$ decreases as $k$ increases.

\subsection{Coefficient Matrix}
In this subsection, we describe the setting of $\bar{\AAA}^{(k)}$.
Each element of matrices used in the iterative linear solver is given by an integer number without fractional bits.
After the initial scaling of the original linear system,  we cast each element of $\AAA$ to an integer number and obtain $\bar{\AAA}_{0}$.
Next, we calculate $\AAA_{1}$ using $\AAA_{1}=\AAA-\bar{\AAA}_{0}$ with FP arithmetic.
Then, we determine a scaling factor $\bar{\alpha}_{1}$ as follows:
\beq
\bar{\alpha}_{1} = \bar{\alpha}_{a} - \lfloor \log_{2}  \max_{i j} | A_{1}(i, j)    |\rfloor.
\eeq{alpha1}
After each element of $\AAA_{1}$ is multiplied by $2^{\bar{\alpha}_{1}}$, it is cast to an intWL number to obtain $\bar{\AAA}_{1}$.
After the same scaling and casting processes are performed repeatedly, the coefficient matrix can be written as
\beq
\AAA = \bar{\AAA}_{0} + \frac{1}{2^{\bar{\alpha}_{1}}} \bar{\AAA}_{1} +\frac{1}{2^{\bar{\alpha}_{2}}} \bar{\AAA}_{2} + \cdots +\frac{1}{2^{\bar{\alpha}_{p}}} \bar{\AAA}_{p}, 
\eeq{A}
because each element of $\AAA$ is an FP number with a finite word length.
Each element of $\bar{\AAA}_{l} (l=0,\ldots,p)$ is an integer number (no fractional bits).
It holds that $\bar{\alpha}_{1} < \bar{\alpha}_{2} < \cdots < \bar{\alpha}_{p}$.
In the $k$-th refinement process, we use a limited number of terms on the right-hand side of (\ref{A}); that is
\beq
\bar{\AAA}^{(k)} =  \bar{\AAA}_{0} + \sum_{l=1}^{s(k)} \frac{1}{2^{\bar{\alpha}_{l}}} \bar{\AAA}_{l},
\eeq{AA}
where $s(k)$ is a parameter for the solver.
When $s(k)=0$, we only use $\bar{\AAA}_{0}$ in the refinement process; that is, $\bar{\AAA}^{(k)} =\bar{\AAA}_{0} $.

\subsection{Refinement Process}
Finally, we introduce an iterative refinement framework for iterative linear solvers that mainly use integer arithmetic, as shown in Fig. \ref{ir}. 
In Fig. \ref{ir}, $\tilde{\xxx}$ is the approximation of $\xxx$ and $S$ is the maximum value of $s(k)$.
We assume that no FP arithmetic is used in the main loop of the iterative solver used in the framework. 
Table \ref{arg} lists the arguments of the iterative linear solver based on integer arithmetic.
The input parameter $d_{f}$ is the number of fractional bits for fixed-point numbers involved in the iterative solver. 
In the program, the input data of the coefficient matrix are represented by integer numbers.
The input of  $\xxx^{(k)}$ is an initial guess for the iterative solver.
The output of $\xxx^{(k)}$ is the (approximate) solution vector of (\ref{solve}), each element of which is an FP number.

\begin{figure}[t]
\begin{center}
\begin{tabular}{|l|}\hline
\\
Initial scaling \\
Calculate $\bar{\AAA}_{0}$, $\bar{\AAA}_{1}$, \ldots, $\bar{\AAA}_{S}$, $\bar{\alpha}_{1}$, $\bar{\alpha}_{2}$, \ldots, $\bar{\alpha}_{S}$\\
for $k=1, 2, \ldots$ \\
\ \  if ($\| \bbb - \AAA \tilde{\xxx}\|  /  \| \bbb \|  <  \epsilon$)  break \\
\ \ Calculate $\bbb'^{(k)}$\\
\ \ Calculate $\gamma^{(k)}$ and $\bbb^{(k)}$ \\
\\
\ \ Integer\_arithmetic\_based\_linear\_solver( arguments )\\
\ \  // to solve $\bar{\AAA}^{(k)} \xxx^{(k)} = \bbb^{(k)}$ \\
\\
\ \ $\tilde{\xxx} \leftarrow \tilde{\xxx} + \gamma^{(k)} \xxx^{(k)}$\\
endfor\\ \hline
\end{tabular}
\caption{Iterative refinement framework using the iterative linear solver based on integer arithmetic}
\label{ir}
\end{center}
\end{figure}

\section{GMRES Solver Using Integer Arithmetic (int-GMRES)}
\subsection{Overview and Data Types of int-GMRES}
In this section, we introduce the GMRES solver based on integer arithmetic that is used in the iterative refinement framework.
We denote the solver by int-GMRES in this paper.

In our solver, each element of the coefficient matrix is given by an integer number (no fractional bits).
The elements of vectors and variables used in the main GMRES iteration loop are given by fixed-point numbers in the Q$d_{m}.d_{f}$ format.  

In the following sections, we use the term "{\it bit shift}." In this paper, left and right shifts with $\beta$ bits refer to multiplication by $2^{\beta}$ and division by $2^{\beta}$, respectively.
These operations for signed integer numbers can be implemented using the {\it shift} operation when the used computer supports a {\it logical shift}. 
In this paper, we assume the use of this type of computer. However, in some computational environments, the result of a shift operation for a signed integer number is "undefined."

Figure \ref{int-gmres} shows the algorithm for the int-GMRES solver of $m$ iterations.
In a practical application, iteration can be terminated when $|\bar{g}_{j+1}|$ is sufficiently small.
In the figure, (FP) represents the statement or calculation based on FP arithmetic, whereas (INT) represents
integer arithmetic.
In the following subsections, we explain for the basic arithmetic of fixed-point numbers and kernels of GMRES, and then present the implementation details.

\subsection{Basic Arithmetic of Fixed-Point Numbers} \label{basic}
In this subsection, we describe the implementation of four basic arithmetics of fixed-point numbers in the Q$d_{m}.d_{f}$ format. 
\subsubsection{Addition and Subtraction}
The addition and subtraction of two fixed point numbers of Q$d_{m}.d_{f}$ are straightforwardly implemented using the integer addition instruction. The obtained integer value directly represents the result in the Q$d_{m}.d_{f}$ format.

\subsubsection{Multiplication}
\begin{figure}[tbp]
\centering
\includegraphics[scale=0.35,clip, bb= 0 75 800 520]{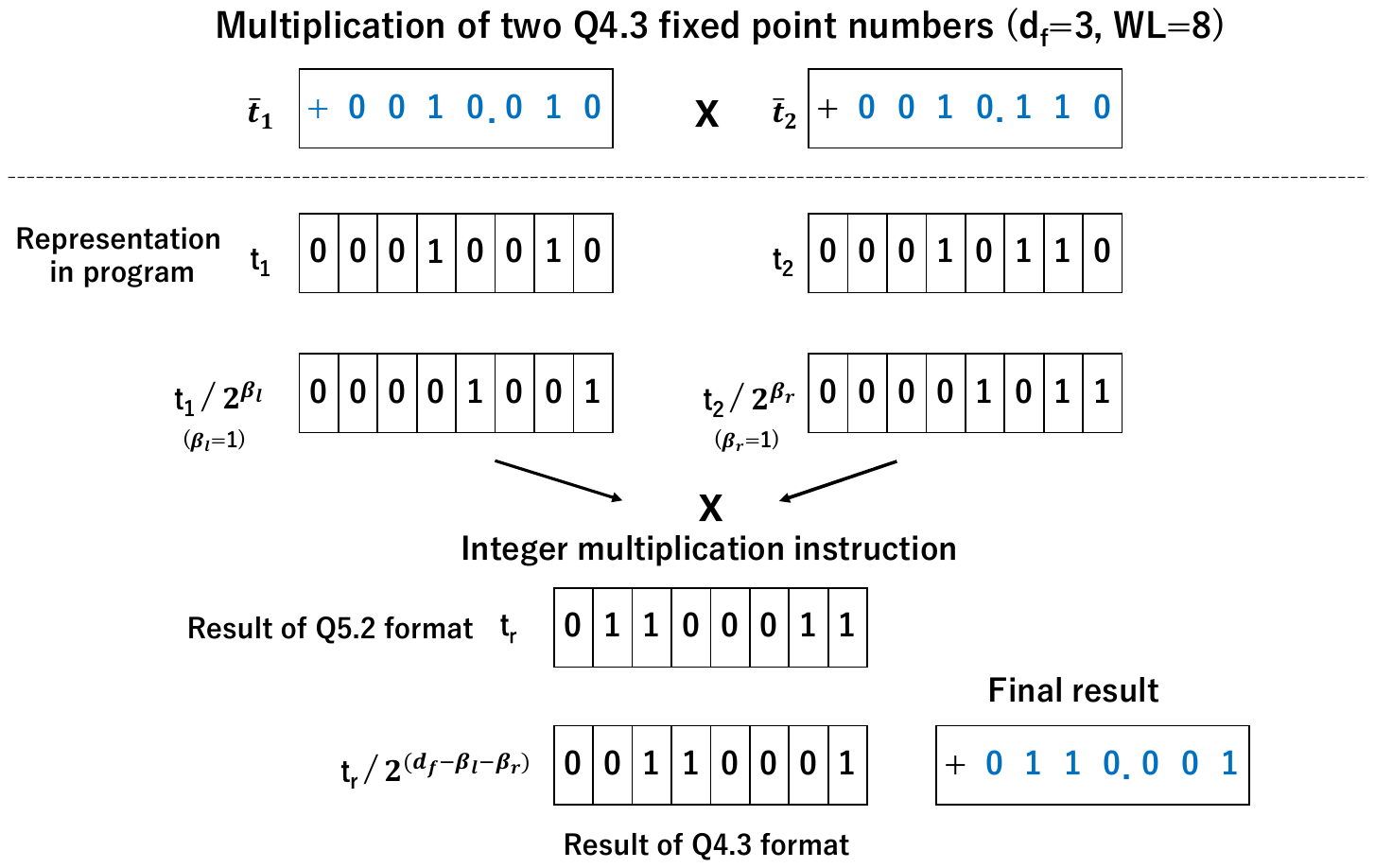}
\caption{Multiplication of fixed-point numbers} 
\label{mul}
\end{figure}

The multiplication of fixed-point numbers is required in various parts of the GMRES program that include calculations of inner products and norms.
Let us consider the multiplication of two fixed-point numbers in the Q$d_{m}.d_{f}$ format: $\bar{t}_{1}$ and $\bar{t}_{2}$. 
We denote the integer representation of $\bar{t}_{1}$ and $\bar{t}_{2}$ in the program by $\sf{t}_{1}$ and $\sf{t}_{2}$, respectively; that is, ${\sf t}_{{\sf 1}}=2^{d_{f}}\cdot \bar{t}_{1}$ and ${\sf t}_{{\sf 2}}=2^{d_{f}}\cdot \bar{t}_{2}$. 
The multiplication procedure for $\bar{t}_{r}=\bar{t}_{1}\bar{t}_{2}$ is given as follows:
After two integer numbers $\sf{t}_{1}$ and $\sf{t}_{2}$ are divided by $2^{\beta_{1}}$ and $2^{\beta_{2}}$, respectively, they are multiplied using the integer instruction. 
The obtained value corresponds to $\bar{t}_{r}$ in the Q$d_{m}'.d_{f}'$ format, where $d_{f}'=2d_{f}-\beta_{1}-\beta_{2}$ and $d_{m}'=WL-d_{f}'-1$.
When we need the result represented as a  Q$d_{m}.d_{f}$ number, the value is divided by $2^{(d_{f}'-d_{f})}$.  
Figure \ref{mul} demonstrates the multiplication procedure of fixed-point numbers.
When the computer supports a {\it logical} shift operation, 
using a C language-like representation, the multiplication in the program, in which the result is represented in the Q$d_{m}.d_{f}$ format, is written as
\beq
{\sf t}_{{\sf r}} = (({\sf t}_{{\sf 1}}>>\beta_{1}) * ({\sf t}_{{\sf 2}}>>\beta_{2}))>>(d_{f}-\beta_{1}-\beta_{2}),  
\eeq{c1c2}
where ${\sf t}_{{\sf r}}$ is the integer representation in the program for $\bar{t}_{r}$.

\begin{table}[t]
\caption{Types of arguments}
\label{arg}
\centering
\normalsize
\begin{tabular}{|c|c|c|}\hline
Arrays, variables & I/O & Number type \\ \hline
$\bar{\AAA}_{0}$, $\ldots$, $\bar{\AAA}_{s}$ & Input& Integer \\ \hline
$\bar{\alpha}_{1}$, $\ldots$, $\bar{\alpha}_{s}$ & Input & Integer \\ \hline
$d_{f}$ & Input & Integer \\ \hline 
$\bbb^{(k)}$ & Input & Floating point \\ \hline
$\xxx^{(k)}$  & Input / Output & Floating point \\ \hline 
\end{tabular}
\end{table}

\begin{figure}[tb]
\begin{center}
\begin{tabular}{|ll|}\hline
1. & Compute $\rrr_{0}=\bbb^{(k)}-\bar{\AAA}^{(k)} \xxx^{(k)}$, \\ 
   &  \ \ \ \ \ \ \ \ \ \ \,  $\vvv_{1}=\rrr_{0}/\|\rrr_{0}\|$ // (FP) \\ 
2.   & Cast $\vvv_{1}$ to $\bar{\vvv}_{1}$ \\
3.   & $\bar{\gggg}=(1, 0, \ldots, 0)^{\top}$ \\
4.  & For $j$=1, 2, \ldots, $m$ \\
5.  & \ \ Compute $\bar{\www}_{j+1}=\bar{\AAA}^{(k)} \bar{\vvv}_{j}$ // (INT) \\
6.  & \ \ For $i=1, \ldots, j$  \\
7. &  \ \ \ \ $\bar{h}_{i, j}=(\bar{\www}_{j+1}, \bar{\vvv}_{i})$ // (INT) \\
8. &  \ \ \ \ $\bar{\www}_{j+1}=\bar{\www}_{j+1} - \bar{h}_{i, j} \bar{\vvv}_{i}$ // (INT)  \\
9. & \ \ Endfor \\
10. &  \ \ $\bar{h}_{j+1,j}=\|\bar{\www}_{j+1}\|$ // (INT)  \\
11. &  \ \ $\bar{\vvv}_{j+1}=\bar{\www}_{j+1} / \bar{h}_{j+1,j}$ // (INT) \\
12. & \ \ For $i=1, \ldots, j-1$  \\
13. & \ \ \ \ $\left( \begin{array}{c}
\bar{h}_{i,j} \\
\bar{h}_{i+1,j}
\end{array} \right)=
\left( \begin{array}{cc}
\bar{c}_{i} & \bar{s}_{i} \\
-\bar{s}_{i} & \bar{c}_{i} \\
\end{array} \right)
\left( \begin{array}{c}
\bar{h}_{i,j} \\
\bar{h}_{i+1,j}
\end{array} \right)
$  \\
   & \ \ \ \ \ \ \ \ \ \ \ \ \ \ \ \ \ \ \ \ \ \ \ \ \ \ \ \ \   // (INT)   \\
14. & \ \ Endfor \\
15. & \ \ $\bar{t}_{mp}=\sqrt{\bar{h}_{j,j}^{2}+\bar{h}_{j+1,j}^{2}}$ // (INT)  \\
16. & \ \ $\bar{c}_{j}=\frac{\bar{h}_{j,j}}{\bar{t}_{mp}}$, $\bar{s}_{j}=\frac{\bar{h}_{j+1,j}}{\bar{t}_{mp}}$  // (INT) \\
17. & \ \ $\bar{g}_{j}=\bar{c}_{j}*\bar{g}_{j}$, \ $\bar{g}_{j+1}=-\bar{s}_{j}*\bar{g}_{j}$ // (INT) \\
18. & \ \ $\bar{h}_{j,j}=\bar{t}_{mp}$ \\
19. & \ \ $\bar{h}_{j+1,j}=0$ // (INT) \\
20. & Endfor \\
21.  & Cast $\bar{\gggg}$ to $\gggg$, and $\bar{\vvv}_{i}$ to $\vvv_{i}$\\
22. & $\yyy= \| \rrr_{0}\| \HHH_{m}^{-1} \gggg$ // (FP) \\
23. & $\xxx^{(k)}=\xxx^{(k)}+\sum_{i=1}^{j}y_{i}\vvv_{i}$ // (FP)\\ \hline
\end{tabular}
\caption{Algorithm for the int-GMRES method ($m$ iterations)}
\label{int-gmres}
\end{center}
\end{figure}

\subsubsection{Division} \label{division}
The division of $\bar{t}_{1}$ by $\bar{t}_{2}$ is implemented as follows:
After the first source operand $\sf{t}_{1}$ is multiplied by $2^{\beta_{1}}$ and the second source operand $\sf{t}_{2}$ is divided by $2^{\beta_{2}}$, the first operand is divided by the second operand using the integer division instruction. 
The resultant variable is multiplied by $2^{(d_{f}-\beta_{1}-\beta_{2})}$ and the final result in the Q$d_{m}.d_{f}$ format is obtained.

\subsubsection{Square Root} \label{label1}
The calculation of a square root is required for the GMRES algorithm. In this subsection, we describe the calculation of the square root of a fixed-point number $\bar{t}_{s}$ in the  Q$d_{m}'.d_{f}'$ format, where $d_{m}'+d_{f}'=WL-1$.
Let $\sf{t}_{s}$ denote the integer representation of $\bar{t}_{s}$ in the program.
We apply the Babylonian square root algorithm for $\sf{t}_{s}$ using integer arithmetic.
The obtained value $\sf{t}_{r} = \sqrt{\sf{t}_{s}}$ is multiplied by $2^{(d_{f}-d_{f}'/2)}$.
The final result provides the integer representation of the square root of $\bar{t}_{s}$ in the Q$d_{m}.d_{f}$ format.

\subsection{Kernels of the GMRES Method}
In this subsection, we describe the implementation of three computational kernels of the GMRES method.

\subsubsection{Inner product} \label{innerproduct}
We consider the inner product of two vectors, each element of which is a fixed-point number in the Q$d_{m}.d_{f}$ format. 
Using the multiplication and addition operations for fixed-point numbers described in Section \ref{basic}, we obtain the result of the inner product as a number in the Q$d_{m}'.d_{f}'$ format.
To obtain a better accuracy in calculations, we typically set $d_{f}'$ to be larger than $d_{f}$.
Therefore, to obtain a result in the Q$d_{m}.d_{f}$ format, the result variable is divided by $2^{d_{f}-d_{f}'}$.
Figure \ref{inner-product} shows a sample code for the inner product. In the figure, \verb|b1| and \verb|b2| correspond to $\beta_{1}$ and $\beta_{2}$ in the procedure for the multiplication, respectively.

\begin{figure}
\begin{verbatim}
cs=0;
for (l=0; l<n; l++){
    cs=cs+(v[l] >> b1)*(w[l] >> b2) ; }
\end{verbatim}
\caption{Calculation of the inner product with the setting $d_{f}'=2d_{f}-\beta_{1}-\beta_{2}$ (only on computers that support a logical shift for a signed integer number) }
\label{inner-product}
\end{figure}

\subsubsection{Norm} \label{norm}
When we calculate a vector norm, we first calculate the inner product of the vector and itself. 
Using the procedure described above, we obtain the result of the inner product in the Q$d_{m}'.d_{f}'$ format.
Then, we calculate its square root using the procedure described in Section \ref{label1}.
Finally, we obtain the norm of the vector which is represented in the Q$d_{m}.d_{f}$ format. 

\subsubsection{Matrix Vector Multiplication} \label{matvec-sec}
Matrix vector multiplication is a main kernel of Krylov subspace methods, in which the GMRES method is classified.
From (\ref{AA}), the kernel consists of $s+1$ matrix vector multiplications:
\beq
\bar{\AAA}^{(k)} \bar{\vvv} =  \bar{\AAA}_{0}\bar{\vvv} + \sum_{l=1}^{s(k)} \frac{1}{2^{\bar{\alpha}_{l}}} \bar{\AAA}_{l}\bar{\vvv},
\eeq{matvec}
where $\bar{\vvv}$ is an $n$-dimensional source vector.
In our implementation, each element of the matrices is given by an integer number, which has no fraction bits. The element of the source and resultant vectors is a fixed-point number in the Q$d_{m}.d_{f}$ format.
Consequently, each matrix vector multiplication $\bar{\AAA}_{l}\bar{\vvv}$ can be performed by a simple integer matrix vector multiplication program.
Each element of $\bar{\AAA}_{l}\bar{\vvv}$ is divided by $2^{\bar{\alpha}_{l}}$, and then added to the corresponding element of the resultant vector.

In the above procedure, it is implied that the result of $\bar{\AAA}_{l}\bar{\vvv}$ does not contribute to the final result when $\bar{\alpha}_{l}$ is substantially  large. Consequently, we estimate that $s$ must be at most 3 or 4 in a practical scenario. When we require more accuracy for the matrix vector multiplication, we should use multiple words for each element of the resultant vector.

\subsection{Implementation Details of int-GMRES and Setting of the Operand Shifts}
In this subsection, we present the details of the int-GMRES solver while paying special attention to setting the parameters in fixed-point number arithmetic.

\subsubsection{Cast of $\vvv_{1}$ to $\bar{\vvv}_{1}$ (l. 2 in Fig. \ref{int-gmres})}
Each element of $\vvv_{1}$ is multiplied by $2^{d_{f}}$ using FP arithmetic.
Then, it is cast to an intWL number.
The obtained integer array that corresponds to $\bar{\vvv}_{1}$ consists of fixed-point numbers in the Q$d_{m}.d_{f}$ format.

\subsubsection{Arnoldi Process (l. 4-11 in Fig. \ref{int-gmres}) }
\paragraph{Line 5 (matrix vector multiplication)}
Line 5 is matrix vector multiplication, which we implement using the method described in Section \ref{matvec-sec}.

\paragraph{Line 7 (inner product)}
Line 7 in Fig. \ref{int-gmres} is the calculation of an inner product, which we implement using the method described in Section \ref{innerproduct}. We suggest a special setting for the operand shift in the multiplication involved in the calculation.
Because $\bar{\vvv}_{i}$ is a normalized vector, the upper $WL-d_{f}-2$ bits of each element of $\bar{\vvv}_{i}$ are always zero.
Considering this feature, we only shift the first source operand which corresponds to $\bar{\www}_{j+1}$; that is $\beta_{2}=0$.

\paragraph{Line 8}
Line 8 involves the multiplication of a vector element by a scalar value and subtraction between two vectors.
Like the inner product in line 7, we only shift the first source operand in the multiplication, considering the profile of $\bar{\vvv}_{i}$.

\paragraph{Line 10 (norm)}
Line 10 is the calculation of the norm, which we implement using the procedure described in Section \ref{norm}. In the multiplication, we naturally set $\beta_{1}=\beta_{2}$.

\paragraph{Line 11}
Line 11 is the division of a vector by a scalar number. We use the procedure for division described in Section \ref{division}.

\subsubsection{Givens Rotation  (l. 12-20 in Fig. \ref{int-gmres}) }
\paragraph{Line 13}
We regard the statement as the inner product of the vectors of two elements.
Therefore, we use the procedure for the inner product.
For the multiplication involved in the procedure, we also use a special setting for the operand shift.
Because the absolute value of $\bar{c}_{i}$ and $\bar{s}_{i}$ is not larger than one, we only shift the second operand corresponding to $\bar{h}_{i,j}$ or  $\bar{h}_{i+1,j}$, that is, $\beta_{1}=0$.

\paragraph{Line 15}
We can implement the statement as the calculation of the norm of the vector of two elements.

\paragraph{Line 16}
We use the procedure for the division of a fixed-point number by another fixed-point number.

\paragraph{Line 17}
Line 17 consists of the multiplication of scalar values. 
The absolute values of $\bar{c}_{i}$ and $\bar{s}_{i}$ are bounded by one, and $g_{j}$ monotonically decreases as the iteration count $j$ increases.
Thus, we do not shift the operands in the multiplication because of the low risk of overflow.

\subsubsection{Update of the Solution Vector  (l. 21-23 in Fig. \ref{int-gmres}) }
\paragraph{Line 21}
We cast each element of the integer arrays for $\bar{\gggg}$ and $\bar{\vvv}_{i}$ to an FP number, which we  then divide by $2^{d_{f}}$. In the practical implementation, we combine these casting operations with the following computations (l. 22-23) to avoid an additional array allocation.

\paragraph{Lines 22 and 23}
We update the output of the int-GMRES solver, that is, $\xxx^{(k)}$ using FP arithmetic.

\subsubsection{Summary of Setting the Parameters}
Table \ref{shift0} summarizes the type of fixed-point numbers, that is, the number of fractional bits, and the quantity of ${\it shift}$ in arithmetic involved in the int-GMRES solver.
In the table, the line number corresponds to the line of the statement in Fig. \ref{int-gmres}, and \#fb represents the number of fractional bits of the fixed-point number used for vectors and variables.

\section{Preconditioning}
Preconditioning is a practically important technique to accelerate the convergence of an iterative solver.
To apply a preconditioning technique to the GMRES solver, we replace two statements (lines 1 and 5) in Fig. \ref{int-gmres} by the following statements:

Line 1'.  \ Compute $\rrr_{0}=\bar{\MMM}^{-1}(\bbb^{(k)}-\bar{\AAA}^{(k)} \xxx^{(k)}),  \\
\hspace{8.95em} \vvv_{1}=\rrr_{0}/\|\rrr_{0}\|$  // (FP)

Line 5'. \ Compute $\bar{\www}_{j+1}=\bar{\MMM}^{-1}\bar{\AAA}^{(k)} \bar{\vvv}_{j}$ // (INT) 

Typically, the preconditioner matrix $\bar{\MMM}$ well approximates the coefficient matrix.
In this paper, we report the application of a standard incomplete LU (ILU), which is precisely ILU(0) preconditioning.
In our solver, the element of the preconditioner matrix is given by an integer number (no fractional bits), which is the same as the coefficient matrix.

\subsection{ILU preconditioning}
In ILU preconditioning, we use the incomplete factorized matrix of the coefficient matrix.
Using FP arithmetic, we incompletely factorize the coefficient matrix as
\beq
\AAA \simeq \LLL \DDD \UUU,
\eeq{ilu}
where $\LLL$ and $\UUU$ have ones for their diagonal elements.
Next, we define two diagonal matrices as follows:
\beq
\DDD_{l}\DDD_{r}=\DDD,
\eeq{dldr}
\beq
D_{l}(i, i)=| d_{ii}|^{1/2},
\eeq{dl}
and
\beq
D_{u}(i, i)=sgn(d_{ii}) | d_{ii}|^{1/2}.
\eeq{dl2}
Then, we introduce two matrices:
\beq
\tilde{\LLL}=\LLL \DDD_{l}
\eeq{l1} 
and
\beq
\tilde{\UUU}=\DDD_{u} \UUU.
\eeq{u1} 
We apply the type cast from float/double to int for each element of $\tilde{\LLL}$ and $\tilde{\UUU}$, and then we obtain  lower and upper triangular matrices, $\bar{\LLL}$ and $\bar{\UUU}$, respectively.
Then, the preconditioner matrix $\bar{\MMM}_{ilu}$ is given by 
\beq
\bar{\MMM}_{ilu}=\bar{\LLL}\bar{\UUU}.
\eeq{ilu2}

The ILU preconditioning step corresponding to lines 1' and 5' is given by forward and backward substitutions.
We can simply use a program for the substitutions in which integer arithmetic is used; that is, if we have a program for substitution based on FP arithmetic, we only change the data type for the matrix and vectors (float/double to int) in the program.

\begin{table*}[t]
\caption{Number of fractional bits of the input and output variables and operand shifts in the calculation}
\label{shift0}
\centering
\begin{tabular}{cccccccc}
 &  & & \multicolumn{4}{c}{Input} & Output \\  \cline{4-8} 
Line \# & Kernel & Arithmetic & \multicolumn{2}{c}{1st (or single) operand} &  \multicolumn{2}{c}{2nd operand} & Result \\ \cline{4-8}
 &  & &  \#fb & Shift & \#fb & Shift & \#fb \\ \hline
Line 5 & Matrix vector  &  Multiplication & 0 & No & $d_{f}$ & No & $d_{f}$ \\ 
         &      Multiplication    &  Addition & $d_{f}$ & No & $d_{f}$ & No & $d_{f}$  \\ \hline
 &  & Multiplication & $d_{f}$ & $/2^{\beta_{1}}$ & $d_{f}$ & No & $2d_{f}-\beta_{1}$ \\
Line 7         &  Inner Product      & Addition & $2d_{f}-\beta_{1}$ & No & $2d_{f}-\beta_{1}$ & No & $2d_{f}-\beta_{1}$  \\
         &                      & Shift & $2d_{f}-\beta_{1}$ & $/2^{d_{f}-\beta_{1}}$  & -  & - & $d_{f}$  \\ \hline
Line 8 & Vector update & Multiplication & $d_{f}$ & $/2^{\beta_{1}}$ & $d_{f}$ & No & $d_{f}$ \\
          &                      & Subtraction & $d_{f}$ & No & $d_{f}$ & No & $d_{f}$  \\ \hline
Line 10 & Norm & Multiplication & $d_{f}$ & $/2^{\beta_{1}}$ & $d_{f}$ & $/2^{\beta_{1}}$ & $2d_{f}-2\beta_{1}$ \\
         &       & Addition & $2d_{f}-2\beta_{1}$ & No & $2d_{f}-2\beta_{1}$ & No & $2d_{f}-2\beta_{1}$  \\
         &        & Square root & $2d_{f}-2\beta_{1}$ & No  & -  & - & $d_{f}$  \\ \hline
Line 11 &   &  Division & $d_{f}$ & $\times2^{\beta_{1}}$ & $d_{f}$ & $/2^{\beta_{2}}$ & $d_{f}$ \\  \hline
 &  & Multiplication & $d_{f}$ & No  & $d_{f}$ & $/2^{\beta_{2}}$ & $2d_{f}-\beta_{2}$ \\
Line 13         &  Inner Product      & Addition & $2d_{f}-\beta_{2}$ & No & $2d_{f}-\beta_{2}$ & No & $2d_{f}-\beta_{2}$  \\
         &                      & Shift & $2d_{f}-\beta_{2}$ & $/2^{d_{f}-\beta_{2}}$  & -  & - & $d_{f}$  \\ \hline
Line 15 &  & Multiplication & $d_{f}$ & $/2^{\beta_{1}}$ & $d_{f}$ & $/2^{\beta_{1}}$ & $2d_{f}-2\beta_{1}$ \\
         &       & Addition & $2d_{f}-2\beta_{1}$ & No & $2d_{f}-2\beta_{1}$ & No & $2d_{f}-2\beta_{1}$  \\
         &        & Square root & $2d_{f}-2\beta_{1}$ & No  & -  & - & $d_{f}$  \\ \hline
Line 16 &   &  Division & $d_{f}$ & $\times2^{\beta_{1}}$ & $d_{f}$ & $/2^{\beta_{2}}$ & $d_{f}$ \\ \hline 
Line 17 &  & Multiplication & $d_{f}$ & No & $d_{f}$ & No & $d_{f}$ \\ \hline
\end{tabular}
\end{table*}

\section{Numerical Result}
\subsection{Computation Environment and Test Problems}
We conducted  numerical tests to evaluate the developed int-GMRES solver.
We evaluated the convergence of the relative residual norm of the solver in comparison with a standard GMRES solver using FP arithmetic.
We performed numerical tests on a node of Fujitsu CX2550 (M4) at the Information Initiative Center, Hokkaido University.
The node was equipped with two Intel 20-core Xeon (Gold 6148) processors  and 384 GB shared memory. 
We wrote the program code in C and used an Intel compiler for the analysis. Logical shift was supported on the computer for a signed integer number.

In the integer arithmetic-based solver, the linear system (\ref{solve}) for the refinement was approximately solved using $m$ iterations of int-GMRES.
For comparison, we also used a standard double precision GMRES($m$) solver.
We set the convergence criterion as the relative residual norm being less than $10^{-8}$.
The relative residual norm was calculated every $m$ iterations using FP arithmetic in both the standard and integer arithmetic-based solvers.
The comparison of the convergence properties of the solvers is performed every $m$ iterations.

For the test problems, we selected ten linear systems from the SuiteSparse Matrix Collection \cite{Florida}.
We selected unsymmetric matrices with various sizes from the collection, for which the standard GMRES solver based on double-precision FP arithmetic worked.
Table \ref{matInfo} lists the properties of selected matrices. The right-hand side vector was given by a vector of ones.

\begin{table}[tbp]
	\centering
	\caption{Matrix information for the test problems}
	\label{matInfo}
	\begin{tabular}{llll}
	\hline
	Data set & Problem type & Dimension & \# nonzero \\
	\hline
atmosmodj & CFD &	1,270,432 & 8,814,880 \\ 
atmosmodl & CFD &	1,489,752	 & 10,319,760\\ 
cage14 & Graph & 1,505,785 & 27,130,349\\ 
CoupCons3D&Structural  problem& 416,800 & 17,277,420\\ 
epb2 & Thermal problem & 25,228 &175,027 \\  
majorbasis & Optimization problem& 160,000 &1,750,416\\ 
memchip	& Circuit simulation & 2,707,524 & 13,343,948\\ 
stomach	& Electro-physical model & 213,360 & 3,021,648 \\ 
torso3 & Finite difference model& 259,156  & 4,429,042 \\ 
wang3 & Semiconductor analysis & 26,064 & 177,168 \\ \hline 
\end{tabular}
\end{table}

\subsection{Results for the Non-Preconditioned GMRES Solver}
The int-GMRES solver based on integer arithmetic requires parameters to be set. The number of fractional bits $d_{f}$ was given by 30. The word length $WL$ was 64, and the 64bit integer (int64) type was used for both fixed-point and integer numbers used in the solver.
The parameter for the coefficient matrix $s$ was given by zero; that is, we only used $\bar{\AAA}_{0}$ in the test.
These settings were also used in the numerical test of the preconditioned GMRES solver.
Table \ref{shift1} lists the setting (number of bits) for operand shifts involved in the calculation of int-GMRES.
Moreover, $\bar{\alpha_{a}}$ was set to 16 for the non-preconditioned solver and 32 for the ILU preconditioned solver.

Table \ref{no-pre-result} shows the number of iterations of a standard GMRES solver using double-precision FP arithmetic and the int-GMRES solver, where "Double" denotes the standard solver.
In the numerical tests, int-GMRES, in addition to the standard solver could solve the problem.
When $m=10$, int-GMRES unexpectedly converged faster in three test cases.
In the {\it wang3} test, which was the worst case for int-GMRES, the solver only required 20\% more iterations than the standard solver. 

\begin{figure}[t]
\centering
\begin{tabular}{c}
\includegraphics[scale=0.53,clip, bb= 50 50 505 340]{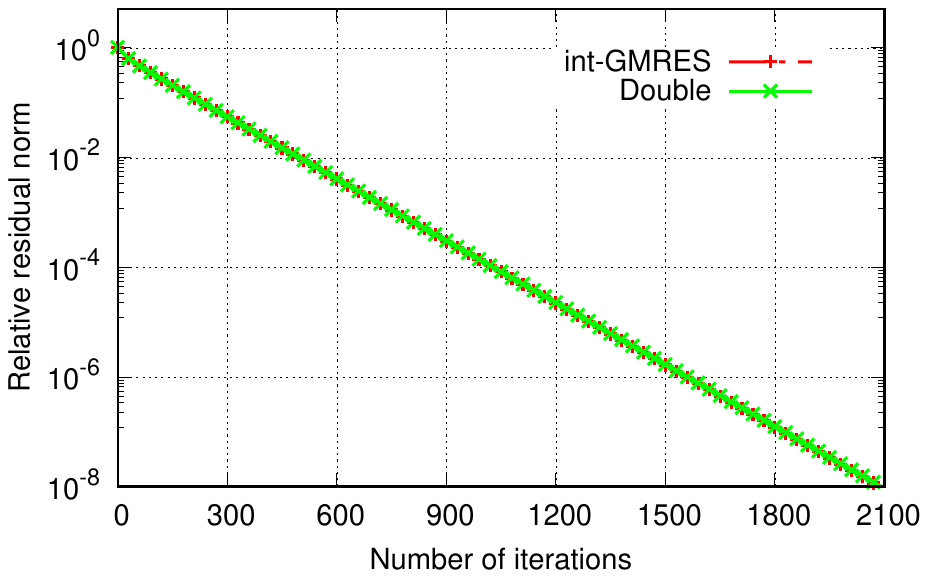} \\
(a) atmosmodj test \\
\\
\includegraphics[scale=0.53,clip, bb= 50 50 505 340]{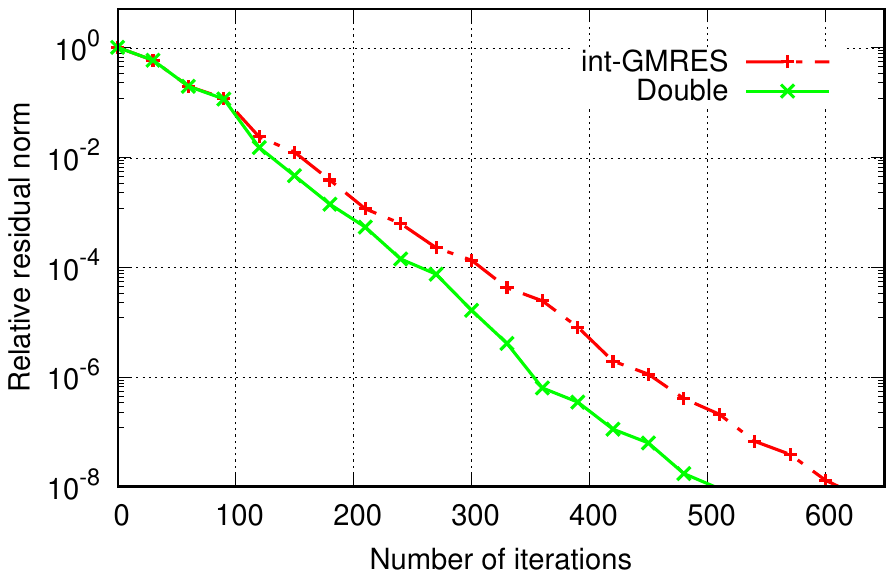} \\
(b) wang3 test \\
\end{tabular}
\caption{Comparison of the convergence behaviors of standard and int-GMRES solvers without preconditioning when $m$ is 30} 
\label{atmos-m30-no}
\end{figure}

When $m=30$, the convergence rates of the standard and integer arithmetic based solvers were comparable for the test cases, except for {\it cage14} and {\it wang3}.
However, int-GMRES only required more iterations for an additional refinement step than the standard solver to converge in {\it cage14}.
In the {\it wang3} test, which was regarded as the worst case for int-GMRES, the solver only required 24\% more iterations than the standard solver. 
Figure \ref{atmos-m30-no} shows the comparison of the convergence rates of the standard double-precision and integer arithmetic-based GMRES solvers.
For the {\it atmosmodj} dataset, the two solvers had an identical convergence rate, which means that the loss of accuracy in int-GMRES did not have a significant influence on solver performance.
In contrast to the result of the {\it atmosmodj} test, the int-GMRES solver had a lower convergence rate than the standard solver in the {\it wang3} test.

\begin{table}[t]
\caption{Number of fractional bits of the input and output variables, and operand shifts for the GMRES solver (no preconditioning)}
\label{shift1}
\centering
\begin{tabular}{ccc}
Line \# & Arithmetic & Setting for operand shift \\ \hline
Line 7  &  Multiplication & $\beta_{1}=16$ \\ \hline
Line 8 & Multiplication & $\beta_{1}=16$  \\ \hline
Line 10 & Multiplication & $\beta_{1}=16$  \\ \hline
Line 11 & Division & $\beta_{1}=16$, $\beta_{2}=14$ \\ \hline
Line 13 & Multiplication & $\beta_{2}=16$ \\ \hline
Line 15 & Multiplication & $\beta_{1}=16$ \\ \hline
Line 16 & Division & $\beta_{1}=16$, $\beta_{2}=14$ \\ \hline
\end{tabular}
\end{table}
\begin{table}[t]
	\centering
	\caption{Number of iterations in no preconditioning case}
	\label{no-pre-result}
	\begin{tabular}{lcccc}
	\hline
     & \multicolumn{2}{c}{$m$=10} & \multicolumn{2}{c}{$m$=30} \\ 
Data set & Double & int-GMRES & Double & int-GMRES  \\
	\hline
atmosmodj & 5,820 & 5,850 & 2,100 & 2,100 \\ 
atmosmodl & 880  & 840 & 420 & 420\\ 
cage14 & 20 & 20 & 30 & 60 \\ 
CoupCons3D & 430 & 430 & 360 & 360 \\ 
epb2 & 820 & 730 & 540 & 540\\  
majorbasis & 90 & 100 & 90 & 90 \\ 
memchip	& 460 & 380 & 300 & 300 \\ 
stomach	& 310 & 310 & 180 & 180 \\ 
torso3 & 150 & 150  & 150  & 150 \\ 
wang3 & 720 & 860 & 510 & 630 \\ \hline

\end{tabular}
\end{table}

\subsection{Results for the Preconditioned GMRES Solver}
When preconditioning is applied to an iterative solver, the convergence rate generally improves.
Thus, we can expect that the residual norm is relatively small and the risk of overflow in the calculation is reduced.
When ILU preconditioning is used, we can avoid the operand shift that sacrifices the accuracy of arithmetics in the solver.
We note that the first source operand shift in the division operation is necessary for improving the calculation accuracy.
Table \ref{shift3} lists the settings of the operand shifts in the preconditioned solver based on integer arithmetic, which is denoted by int-ILU-GMRES.

\begin{figure}[t]
\centering
\includegraphics[scale=0.53,clip, bb= 50 50 505 340]{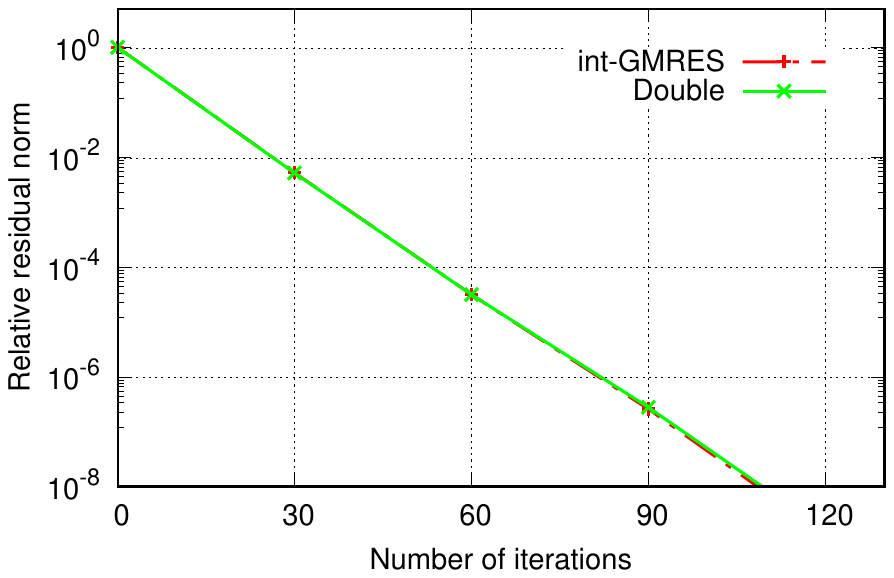}
\caption{Comparison of the convergence behaviors of standard and int-GMRES solvers with ILU preconditioning in the wang3 test when $m$ is 30} 
\label{wang-m30-ilu}
\end{figure}

Table \ref{ilu-pre-result} shows the number of iterations of the ILU-GMRES solver using double-precision FP arithmetic and the int-ILU-GMRES solver.
When compared with the non-preconditioned solver, both solvers attained significant improvement in convergence.
Moreover, the convergence rates of the two solvers were comparable.
The int-ILU-GMRES only required more iterations for an additional refinement step than the standard solver in some test cases.
Figure \ref{wang-m30-ilu} shows that both solvers had identical convergence behavior of the relative residual norm  in the {\it wang3} test.

\subsection{Discussions}
\subsubsection{Preconditioning}
Preconditioning is important in the context of iterative solvers based on integer arithmetic, because it reduces the risk of overflow.
Consequently, we can decrease the number of bits of the operand shift, which improves the accuracy of arithmetic.
For the non-preconditioned solver, we investigated an auto-tuning technique for the shift.
However, it proved to be hardly necessary in the preconditioning case.
In the implementation of int-ILU-GMRES, we could avoid using the operand shift which sacrificed the accuracy. A similar effect was also confirmed in the Gauss--Seidel preconditioning case. 
(Because of the page limit, the numerical result is not shown in this paper.)

\subsubsection{Condition of the Problems}
Because we selected test problems (matrices) for which a non-preconditioned GMRES solver using FP arithmetic attained convergence, the problems were not heavily ill-conditioned.
Consequently, the int-GMRES solver also solved the problems.
It is possible that problems exist that the standard FP arithmetic solver can solve but int-GMRES cannot.
However, as far as we have tested,  it seems not to be an easy task to seek such a problem; that is, the int-GMRES solver used with the iterative refinement technique may have comparable solver performance to the standard FP arithmetic solver.

\begin{table}[t]
\caption{Number of fractional bits of the input and output variables and operand shifts for the ILU-GMRES solver}
\label{shift3}
\centering
\begin{tabular}{ccc}
Line \# & Arithmetic & Setting for operand shift \\ \hline
Line 7  &  Multiplication & $\beta_{1}=0$ \\ \hline
Line 8 & Multiplication & $\beta_{1}=0$  \\ \hline
Line 10 & Multiplication & $\beta_{1}=0$  \\ \hline
Line 11 & Division & $\beta_{1}=30$, $\beta_{r}=0$ \\ \hline
Line 13 & Multiplication & $\beta_{2}=0$ \\ \hline
Line 15 & Multiplication & $\beta_{1}=0$ \\ \hline
Line 16 & Division & $\beta_{1}=30$, $\beta_{2}=0$ \\ \hline
\end{tabular}
\end{table}

\begin{table}[tbp]
	\centering
	\caption{Number of iterations in the ILU preconditioning case}
	\label{ilu-pre-result}
	\begin{tabular}{lcccc}
	\hline
     & \multicolumn{2}{c}{$m$=10} & \multicolumn{2}{c}{$m$=30} \\ 
Data set & Double & int-GMRES & Double & int-GMRES  \\
	\hline
atmosmodj & 610 & 610 & 300 & 300 \\ 
atmosmodl & 140  & 140 & 120 & 120\\ 
cage14 & 10 & 20 & 30 & 60 \\ 
CoupCons3D & 140 & 150 & 30 & 60 \\ 
epb2 & 50 & 50 & 60 & 60\\  
majorbasis & 20 & 20 & 30 & 60 \\  
stomach	& 20 & 20 & 30 & 60 \\ 
torso3 & 40 & 40  & 30 & 60\\ 
wang3 & 180 & 180 & 120 & 120 \\ \hline

\end{tabular}
\end{table}

\section{Related Works}
In this section, we introduce several papers that discuss mixed-precision linear solvers using the iterative refinement technique.
The survey paper \cite{Goddeke} by D. G\"{o}ddeke et al. provides a good introduction to the mixed-precision iterative refinement algorithm framework. The paper \cite{Anzt} by Anzt et al. is another early work on a mixed-precision linear solver, in which the authors reported a GPU implementation of an error correction solver using the GMRES method and showed the effectiveness of their approach in CFD applications. A. Haidar et al. reported the development of an architecture-specific algorithm and highly tuned implementations for the latest GPUs of mixed-precision iterative refinement solvers in \cite{Haidar}. Their solver that involved LU factorization was targeted at a linear system with a dense coefficient matrix. Carson et al. presented a general algorithm for iterative refinement with three precisions and its error analysis in \cite{Carson}. Moreover, the Exascale Computing Project Multiprecision Effort Team (Lead: Hartwig Anzt) recently opened its technical report to the public, which provides a comprehensive review of mixed-precision computing \cite{ECP}.  

Next, we briefly mention analyses based on integer arithmetic (fixed-point numbers). 
Currently, integer arithmetic is often used in machine learning and artificial intelligence applications. LU factorization based on integer arithmetic for these applications is given in \cite{ECP}. Numerical linear algebra algorithms based on fixed-point numbers have also been investigated in the context of signal processing \cite{Niko, Pradhan}. The difference between the present research and these papers is in the investigation and development of the GMRES method based on integer arithmetic. 

\section{Conclusions}
In this paper, we developed a GMRES solver based on integer arithmetic, denoted by int-GMRES. The int-GMRES solver was used with an iterative refinement technique to attain a solution as accurate as that of a normal linear solver based on FP arithmetic. We also developed an ILU preconditioned int-GMRES solver. In integer arithmetic (fixed-point number) computing, it is important to avoid overflow in calculations. We explained how the operands are adjusted (logically shifted) in the calculation considering the characteristics of the GMRES method.
We conducted numerical tests using matrices from SuiteSparse Matrix Collections. The numerical results demonstrated that the int-GMRES solver had comparable solver performance in terms of convergence to the standard FP solver. Moreover, we found that preconditioning was important for the solver using integer arithmetic to avoid overflow.

In the future, we will evaluate solver performance in terms of timing on the model of new computing devices in which integer arithmetic has advantages over conventional computing devices for calculation speed or power consumption.


\begin{thebibliography}{1}


\bibitem{Vetter} J. S. Vetter, E. P. DeBenedictis and T. M. Conte, ``Architectures for the Post-Moore Era," IEEE Micro, vol. 37, pp. 6--8, 2017.

\bibitem{Shalf} J. Shalf, ``The future of computing beyond Moore's law," Phil. Trans. R. Soc. A. , 378 (2166), 20190061, 2020.


\bibitem{Sato} R. Sato, Y. Hatanaka, Y. Ando, M. Tanaka, A. Fujimaki, K. Takagi, N. Takagi, 
``High-speed operation of random-access-memory-embedded microprocessor with minimal instruction set architecture based on rapid single-flux-quantum logic," IEEE Trans. Appl. Supercond., vol. 27, pp. 1--5, 2017.

\bibitem{ishida}
K. Ishida, M. Tanaka, T. Ono, K. Inoue, ``Towards ultra-high-speed cryogenic single-flux-quantum computing,"  IEICE Trans. Electron., vol. E101-C, pp. 359--369, 2018.

\bibitem{Saad} Y. Saad, Iterative Methods for Sparse Linear Systems, (Second ed.). Philadelphia, PA, SIAM, 2003. 

\bibitem{Goddeke} D. G\"{o}ddeke, R. Strzodka, and S. Turek, ``Performance and accuracy of hardware-oriented native-, emulated- and mixed-precision solvers in FEM simulations," Int. J. Parallel Emergent Distrib. Syst., vol. 22, pp. 221--256, 2007.


\bibitem{Florida} T. A. Davis, and Y. Hu, ``The university of Florida sparse matrix collection," ACM Trans. Math. Software. 38, pp. 1--25, 2011.


\bibitem{Anzt} H. Anzt, V. Heuveline, and B. Rocker, ``An error correction solver for linear systems: evaluation of mixed precision implementations," VECPAR 2010,  LNCS, vol. 6449, Springer, 2011. 

\bibitem{Haidar} 
A. Haidar, S. Tomov, J. Dongarra, and N. J. Higham, ``Harnessing GPU tensor cores for fast FP16 arithmetic to speed up mixed-precision iterative refinement solvers," Proc. SC18,  Intl. Conf. High Performance Comput., Networking, Storage and Analysis, pp. 603--613, 2018.

\bibitem{Carson}
E. Carson, and N. J. Higham, ``Accelerating the solution of linear systems by iterative refinement in three precisions," SIAM J. Sci. Comput., vol. 40, pp. A817--A847, 2018.

\bibitem{ECP} A. Abdelfattah, H. Anzt, E.G. Boman, E. Carson, T. Cojean, J. Dongarra, M. Gates, T. Gr\"{u}tzmacher, N.J. Higham, S. Li, N. Lindquist, Y. Liu, J. Loe, P. Luszczek, P. Nayak, S. Pranesh, S. Rajamanickam, T. Ribizel, B. Smith, K. Swirydowicz, S. Thomas, S. Tomov, Y.M. Tsai, I. Yamazaki, U.M. Yang, ``A Survey of Numerical Methods Utilizing Mixed Precision Arithmetic," arXiv preprint arXiv:2007.06674, 2020.


\bibitem{Niko}
Z. Nikoli\'{c}, H. T. Nguyen, and G. Frantz, ``Design and implementation of numerical linear algebra algorithms on fixed point DSPs,'' EURASIP Journal on Advances in Signal Processing, 087046, 2007.

\bibitem{Pradhan} T. Pradhan, B. Kabi, and A. Routray, ``Fixed-point Hestenes algorithm for singular value decomposition of symmetric matrices," Proc. 2013 Intl. Conf. Electronics, Signal Processing and Communication Systems, 2013.

\end{thebibliography}
\end{document}